\documentclass[aoas,preprint]{imsart}

\RequirePackage[OT1]{fontenc}
\RequirePackage{amsthm,amsmath}
\RequirePackage[numbers]{natbib}
\RequirePackage[colorlinks,citecolor=blue,urlcolor=blue]{hyperref}
\usepackage[pdftex]{graphicx}

\startlocaldefs
\numberwithin{equation}{section}
\theoremstyle{plain}

\endlocaldefs

\begin{document}

\begin{frontmatter}
\title{Exponential two-armed bandit problem} \runtitle{Exponential two-armed bandit}

\begin{aug}
\author{\fnms{Alexander} \snm{Kolnogorov}\thanksref{m1}\ead[label=e1]{Alexander.Kolnogorov@novsu.ru}},
\and
\author{\fnms{Denis} \snm{Grunev}\thanksref{m1}
\ead[label=e2]{   leonard\_max@mail.ru } }

\runauthor{A. Kolnogorov and D. Grunev}

\affiliation{Yaroslav-the-Wise Novgorod State
University\thanksmark{m1}}

\address{41 B.Saint-Petersburgskaya Str.,
Velikiy Novgorod, Russia, 173003\\
Applied Mathematics and Information Science Department\\
\printead{e1}\\
\printead{e2}}

\end{aug}

\begin{abstract}
We consider exponential two-armed bandit problem in which incomes
are described by exponential distribution densities. We develop
Bayesian approach and present recursive equation for determination
of Bayesian strategy and Bayesian risk. In the limiting case as
the control horizon goes to infinity, we obtain the second order
partial differential equation in the domain of ``close
distributions''. Results are compared with Gaussian two-armed
bandit. It turned out that exponential and Gaussian two-armed
bandits have the same description in the limiting case. Since
Gaussian two-armed bandit describes the batch processing, this
means that in case of exponential two-armed bandit batch
processing does not enlarge Bayesian risk in comparison with
one-by-one optimal processing as the total number of processed
data items goes to infinity.
\end{abstract}



\end{frontmatter}

\def \mE{\mathrm{E}}
\def \mD{\mathrm{Var}}
\def \tR{\tilde{R}}
\def \tr{\tilde{r}}
\def \tx{\hat{x}}
\def \eps{\varepsilon}

\section{Introduction}\label{intro}

We consider the two-armed bandit problem (see, e.g. \cite{BF, PS})
in the following setting. Let $\xi_n$, $n=1,2,\dots,N$, be a
controlled random process which values are interpreted as incomes,
depend only on currently chosen actions $y_n \in \{1,2\}$ and are
described by exponential distribution density
\begin{gather}\label{a1}
f(x|m_\ell)= \left\{
\begin{array}{l}
m_\ell^{-1}\exp(-x m_\ell^{-1}), \ \ x\geq 0,\\
0, \qquad \qquad \qquad \qquad x<0,
\end{array}
\right.
\end{gather}
if $y_n=\ell$, $\ell=1,2$. Here $m_1, m_2$ are one-step
mathematical expectations of income and a vector parameter
$\theta=(m_1,m_2)$ completely describes exponential two-armed
bandit. We assume that the set of admissible values of parameters
$\Theta=\{\theta\}$ is a priori known.

A control strategy $\sigma$ generally assigns a random choice of
the action at the point of time $n+1$ depending on currently
observed history of the process. For exponential distributions of
incomes the history consists of cumulative numbers $n_1, n_2$
($n_1+n_2=n$) of both actions applications and corresponding
cumulative incomes $X_1, X_2$. If one knew both $m_1,m_2$, he
should always choose the action corresponding to the largest of
them, his total expected income on the control horizon $N$ would
thus be equal to $N\max(m_1,m_2)$. But if he uses strategy
$\sigma$, his total expected income is less than maximal by the
value
\begin{gather}\label{a2}
L_N(\sigma,\theta)=N
\max(m_1,m_2)-\mE_{\sigma,\theta}\left(\sum_{n=1}^N \xi_n \right)
\end{gather}
which is called the regret. Here $\mE_{\sigma,\theta}$ denotes the
mathematical expectation with respect to the measure generated by
strategy $\sigma$ and parameter $\theta$.

Let's assign a prior distribution density
$\mu(\theta)=\mu(m_1,m_2)$ on the set of parameters $\Theta$.
Corresponding Bayesian risk is defined as follows
\begin{gather}\label{a3}
R_N(\mu)= \inf_{\{\sigma\}}\int_\Theta
L_N(\sigma,\theta)\mu(\theta) d\theta,
\end{gather}
the optimal strategy $\sigma^B$ is called Bayesian strategy. Note
that Bayesian approach allows to determine Bayesian risk and
Bayesian strategy by solving recursive Bellman-type equation for
arbitrary prior distribution. The minimax risk on the set $\Theta$
is defined as
\begin{gather}\label{a4}
R^M_N(\Theta)= \inf_{\{\sigma\}}\sup_\Theta L_N(\sigma,\theta),
\end{gather}
corresponding optimal strategy $\sigma^M$ is called minimax
strategy. There is no a direct method of determining minimax
strategy and minimax risk. However, one can determine them using
the main theorem of the theory of games according to which the
following equality holds
\begin{gather}\label{a5}
R^M_N(\Theta)= R_N(\mu_0)= \sup_{\{\mu\}}R_N(\mu),
\end{gather}
i.e. minimax risk is equal to the Bayesian calculated with respect
to the worst-case prior distribution and minimax strategy is equal
to corresponding Bayesian strategy.

There are some different approaches to the two-armed bandit
problem. We refer here to \cite{Sragovich}, \cite{CBL},
\cite{Koln18} and references therein.

The rest of the paper is the following. Recursive Bellman-type
equation for determining  Bayesian risk and Bayesian strategy is
presented in Section~\ref{equ1}. Another version of recursive
equation is presented in Section~\ref{equ2}. In Section~\ref{lim},
we obtain a second order partial differential equation in the
limiting case as $N \to \infty$. In Section~\ref{Gaussian} we
compare exponential and Gaussian two-armed bandits. It turned out
that they have the same description in the limiting case. Since
Gaussian two-armed bandit describes the batch processing (see,
e.g.~\cite{Koln18}), this means that in case of exponential
two-armed bandit batch processing does not enlarge Bayesian risk
in comparison with one-by-one optimal processing asymptotically as
$N \to \infty$.

\section{Recursive equation}\label{equ1}

Let's consider control strategies
$\{\sigma_\ell(X_1,t_1,X_2,t_2)\}$ which are defined by a
condition
\begin{gather*}
\Pr(y_{n+1}=\ell|X_1,n_1,X_2,n_2)=\sigma_\ell(X_1,n_1,X_2,n_2),
\end{gather*}
where $n_1, n_2$ are current cumulative times of both actions
applications, $X_1, X_2$ are corresponding current cumulative
incomes. The posterior distribution at the point of time
$n=n_1+n_2$ is calculated as
\begin{gather}\label{b1}
\mu(m_1,m_2|X_1,n_1,X_2,n_2)=\frac{f(X_1,n_1|m_1)f(X_2,n_2|m_2)
\mu(m_1,m_2)}{\mu(X_1,n_1,X_2,n_2)},
\end{gather}
where
\begin{gather}\label{b2}
\mu(X_1,t_1,X_2,t_2)\\ =\iint_\Theta
f(X_1,n_1|m_1)f(X_2,n_2|m_2)\mu(m_1,m_2)dm_1 dm_2.\nonumber
\end{gather}
Here $f(X,n|m)$ is defined as
\begin{gather*}
f(X,n|m)= \left\{
\begin{array}{l}
\displaystyle{\frac{X^{n-1}}{m^n (n-1)!}}  \exp(-X m^{-1}), \ \ X \geq 0,\\
0, \qquad \qquad \qquad \qquad \qquad \quad X<0,
\end{array}
\right.
\end{gather*}
Note that $f(X,1|m)=f(X|m)$ is exponential distribution density.
Let's put $f(0,0|m)=1$. Then~\eqref{b1} remains correct if $n_1=0$
and/or $n_2=0$. Denote $x^+=\max(x,0)$. Using \eqref{a1},
\eqref{a2} we obtain the following recursive Bellman-type equation
for determining Bayesian risk \eqref{a3} with respect to the
posterior distribution \eqref{b1}:
\begin{gather}\label{b3}
R(X_1,n_1,X_2,n_2)= \min (R^{(1)}(X_1,n_1,X_2,n_2),
R^{(2)}(X_1,n_1,X_2,n_2)),
\end{gather}
where
\begin{gather}\label{b4}
R^{(1)}(X_1,n_1,X_2,n_2)= R^{(2)}(X_1,n_1,X_2,n_2)=0
\end{gather}
if $n_1+n_2=N$ and then
\begin{gather}\nonumber
R^{(1)}(X_1,n_1,X_2,n_2)\\
\nonumber =\displaystyle{\iint_\Theta}\mu(m_1,m_2|X_1,n_1,X_2,n_2)
\times \Big( (m_2-m_1)^+ \\
\nonumber +\displaystyle{\int_{0}^\infty} R(X_1+Y,n_1+1,X_2,n_2)
f(Y|m_1) dY \Big)
dm_1 dm_2,\\
\label{b5}
R^{(2)}(X_1,n_1,X_2,n_2)\\
\nonumber =\displaystyle{\iint_\Theta}\mu(m_1,m_2|X_1,n_1,X_2,n_2)
\times \Big( (m_1-m_2)^+\\
+\displaystyle{ \int_{0}^\infty } R(X_1,n_1,X_2+Y,n_2+1) f(Y|m_2)
dY \Big) dm_1 dm_2. \nonumber
\end{gather}
Here  $\{R^{(\ell)}(X_1,n_1,X_2,n_2)\}$ denote expected losses if
initially the $\ell$-th action is applied at the point of time
$n+1$ and then control is optimally implemented ($\ell=1,2$).
Bayesian risk \eqref{a3} is as follows
\begin{gather}\label{b6}
   R_N(\mu)=R(0,0,0,0).
\end{gather}

Equation \eqref{b3}--\eqref{b5} allow to determine Bayesian
strategy, too. Bayesian strategy prescribes to choose $\ell$-th
action if $R^{(\ell)}(X_1,n_1,X_2,n_2)$ has smaller value. In case
of a draw $R^{(1)}(X_1,n_1,X_2,n_2)=R^{(2)}(X_1,n_1,X_2,n_2)$ the
choice of the action may be arbitrary.

Given $N$ large enough, consider the strategy which at the start
of the control $2 n_0$ times equally applies both actions and then
optimally controls. In this case
\begin{gather}\label{b7}
   R_N(\mu)=n_0 \iint_\Theta |m_2-m_1| \mu(m_1,m_2) dm_1 dm_2\\
   + \iint_0^\infty R(n_0,X_1,n_0,X_2) \mu(X_1,n_1,X_2,n_2) dX_1
   dX_2.\nonumber
\end{gather}

\section{One more version of recursive equation}\label{equ2}

In this section, we obtain another version of recursive
Bellman-type equation. Let's denote
\begin{gather*}
   \tR(X_1,n_1,X_2,n_2)=R(X_1,n_1,X_2,n_2)\times
   \mu(X_1,n_1,X_2,n_2),
\end{gather*}
where $\{R(X_1,n_1,X_2,n_2)\}$ are Bayesian risks calculated with
respect to the posterior distribution~\eqref{b1} and
$\{\mu(X_1,n_1,X_2,n_2)\}$ are defined in~\eqref{b2}. Then the
following recursive equation holds
\begin{gather}\label{c1}
\tR(X_1,n_1,X_2,n_2)= \min (\tR^{(1)}(X_1,n_1,X_2,n_2),
\tR^{(2)}(X_1,n_1,X_2,n_2)),
\end{gather}
where
\begin{gather}\label{c2}
\tR^{(1)}(X_1,n_1,X_2,n_2)= \tR^{(2)}(X_1,n_1,X_2,n_2)=0
\end{gather}
if $n_1+n_2=N$ and then
\begin{gather}\nonumber
\tR^{(1)}(X_1,n_1,X_2,n_2)=G^{(1)}(X_1,n_1,X_2,n_2)\\
\label{c3}+\displaystyle{ \int_{0}^\infty}
\tR(X_1+Y,n_1+1,X_2,n_2)\times \frac{n_1 X_1^{n_1-1}}
{(X_1+Y)^{n_1}} dY,\\
\nonumber \tR^{(2)}(X_1,n_1,X_2,n_2)=G^{(2)}(X_1,n_1,X_2,n_2)\\
+\displaystyle{\int_{0}^\infty} \tR(X_1,n_1,X_2+Y,n_2+1) \times
\frac{n_2 X_2^{n_2-1}} {(X_2+Y)^{n_2}}dY. \nonumber
\end{gather}
Here
\begin{gather}\nonumber
G^{(1)}(X_1,n_1,X_2,n_2)\\= \label{c4} \iint_\Theta(m_2-m_1)^+
f(X_1,n_1|m_1)
f(X_2,n_2|m_2) \mu(m_1,m_2) dm_1 dm_2,\\
\nonumber G^{(2)}(X_1,n_1,X_2,n_2)\\=\iint_\Theta(m_1-m_2)^+
f(X_1,n_1|m_1) f(X_2,n_2|m_2) \mu(m_1,m_2) dm_1 dm_2.\nonumber
\end{gather}
Bayesian strategy prescribes to choose  $\ell$-th action if
$\tR^{(\ell)}(X_1,n_1,X_2,n_2)$ has smaller value. In case of a
draw the choice of the action is arbitrary. Bayesian
risk~\eqref{a3} is calculated by the formula
\begin{gather}\label{c5}
   R_N(\mu)=\tR(0,0,0,0).
\end{gather}

Given $N$ large enough, consider the strategy which at the start
of the control $2 n_0$ times equally applies both actions and then
optimally controls. In this case
\begin{gather}\label{c6}
   R_N(\mu)=n_0 \iint_\Theta |m_2-m_1| \mu(m_1,m_2) dm_1 dm_2\\
   + \iint_0^\infty \tR(n_0,X_1,n_0,X_2) dX_1 dX_2.\nonumber
\end{gather}

Formulas~\eqref{c1}--\eqref{c6} follow
from~\eqref{b3}--\eqref{b7}.

\section{A limiting description}\label{lim}

In this section, we present a limiting description by the second
order partial differential equation.  We consider the domain of
``close distributions'', satisfying condition $|m_1-m_2| \le c
N^{-1/2}$ with $c$ large enough but independent from $N$, because
just in this domain the maximum expected losses take place. Denote
$\eps=N^{-1}$, $\delta=N^{-1/2}$, so that $\delta^2=\eps$. Note
that one-step expected income and variance of exponential
two-armed bandit are the following
\begin{gather*}
    \mE(\xi_n||y_n=\ell)=m_\ell, \qquad \mD
    (\xi_n||y_n=\ell)=D_\ell=m^2_\ell,
\end{gather*}
$\ell=1,2$. In the domain of distributions such that $m_1, m_2$
are close to $m$, let's put
\begin{gather}\nonumber
m_\ell=m+(D/N)^{1/2} v_\ell,\\ \nonumber
X_\ell=n_\ell m+x_\ell (DN)^{1/2},\\
\label{d1}Y= m+y \delta (DN)^{1/2},\\ \nonumber
n_\ell=t_\ell N,\\
\nonumber  \mu(m_1,m_2)=(D/N)^{-1}\rho(v_1,v_2),\\ \nonumber
\tx_\ell=x_\ell/t_\ell,\nonumber
\end{gather}
where $D=m^2$, $|x_\ell| \ll N^{1/2}$, $|y| \ll N^{1/2}$; $\ell=1,2$. 
Let's estimate functions in~\eqref{c4}. If $n_\ell$ is large
enough then according to central limit theorem we have
\begin{gather*}
    f(X_\ell,n_\ell|m_\ell)= \frac{1} {(2 \pi D n_\ell)^{1/2}}
    \exp\left\{- \frac{(X_\ell - m_\ell n_\ell)^2}{2 D n_\ell} \right\}(1+o(1))\\
    =\frac{\delta} {(2 \pi D t_\ell)^{1/2}}
    \exp\left\{- \frac{(x_\ell - v_\ell t_\ell)^2}{2 t_\ell}
    \right\}(1+o(1)).
\end{gather*}
Hence, for functions in~\eqref{c4} one derives
\begin{gather}\label{d2}
G^{(\ell)}(X_1,n_1,X_2,n_2)=\eps \left(DN\right)^{-1/2}
g^{(\ell)}(x_1,t_1,x_2,t_2)(1+o(1)),
\end{gather}
where
\begin{gather}\nonumber
g^{(1)}(x_1,t_1,x_2,t_2)\\
\nonumber =\iint_\Theta \frac{(v_2-v_1)^+} {2 \pi (t_1 t_2)^{1/2}}
    \exp\left\{- \frac{(x_1 - v_1 t_1)^2}{2 t_1} - \frac{(x_2 - v_2 t_2)^2}{2 t_2}
    \right\} \rho(v_1,v_2) dv_1 dv_2,\\
\label{d3} g^{(2)}(x_1,t_1,x_2,t_2)\\=\iint_\Theta
\frac{(v_1-v_2)^+} {2 \pi (t_1 t_2)^{1/2}}
    \exp\left\{- \frac{(x_1 - v_1 t_1)^2}{2 t_1} - \frac{(x_2 - v_2 t_2)^2}{2 t_2}
    \right\} \rho(v_1,v_2) dv_1 dv_2. \nonumber
\end{gather}

Let's estimate factors in integrals of~\eqref{c3}. First, one
derives
\begin{gather}\label{d4}
\frac{n_\ell X_\ell^{n_\ell-1}} {(X_\ell+Y)^{n_\ell}}=
\frac{n_\ell}{X_\ell}\times \left(
    1+\frac{Y}{X_\ell}\right)^{-n_\ell}=m^{-1} f_\ell(1+y)
\end{gather}
with
\begin{gather*}
f_\ell(1+y)=\frac{1}{1+\delta \tx_\ell}
    \times \left(1+\frac{1+y}{N t_\ell(1+\delta \tx_\ell)}  \right)^{-t_\ell
    N}.
\end{gather*}
Let's now estimate $f_\ell(1+y)$. Since
\begin{gather*}
\log \left(1+\frac{1+y}{Nt_\ell(1+\delta \tx_\ell)}
\right)^{-t_\ell N}\\ =-t_\ell N \left(\frac{1+y}{Nt_\ell(1+\delta
\tx_\ell)}-
\frac{(1+y)^2}{2 N^2 t_\ell^2(1+\delta \tx_\ell)^2} \right)+o(\eps)\\
=-\frac{1+y}{1+\delta \tx_\ell}+ \frac{\eps (1+y)^2}{2 t_\ell
(1+\delta \tx_\ell)^2}+o(\eps),
\end{gather*}
and $\eps (1+y)^2  \ll 1$, we obtain that
\begin{gather*}
f_\ell(1+y)=\frac{1}{1+\delta
\tx_\ell}\exp\left(-\frac{1+y}{1+\delta
\tx_\ell}\right)\left(1+\frac{\eps (1+y)^2}{2 t_\ell (1+\delta
\tx_\ell)^2}+o(\eps)\right).
\end{gather*}
So, one can verify that
\begin{gather}\nonumber
\int_{-1}^\infty f_\ell(1+y) dy=\int_{0}^\infty f_\ell(z)
dz=1+\eps/t_\ell+o(\eps),\\ \label{d5}
\int_{-1}^\infty y
f_\ell(1+y) dy=\int_{0}^\infty (z-1) f_\ell(z)
dz=\delta \tx_\ell+o(\delta),\\
\int_{-1}^\infty y^2 f_\ell(1+y) dy=\int_{0}^\infty (z-1)^2
f_\ell(z) dz=1+o(1). \nonumber
\end{gather}

Let's put
\begin{gather}\label{d6}
    \tR(X_1,n_1,X_2,n_2)=\left(DN\right)^{-1/2} \tr(x_1,t_1,x_2,t_2)
\end{gather}

Using \eqref{d1}--\eqref{d4} and \eqref{d6}, one derives
from~\eqref{c1}--\eqref{c4} the integro-difference equation
\begin{gather}\label{d7}
\tr(x_1,t_1,x_2,t_2)= \min (\tr^{(1)}(x_1,t_1,x_2,t_2),
\tr^{(2)}(x_1,t_1,x_2,t_2)),
\end{gather}
where
\begin{gather}\label{d8}
\tr^{(1)}(x_1,t_1,x_2,t_2)= \tr^{(2)}(x_1,t_1,x_2,t_2)=0
\end{gather}
if $t_1+t_2=1$ and then
\begin{gather} \nonumber
\tr^{(1)}(x_1,t_1,x_2,t_2)=\eps g^{(1)}(x_1,t_1,x_2,t_2)\\
\nonumber +\displaystyle{ \int_{-1}^\infty} \tr(x_1+\delta
y,t_1+\eps,x_2,t_2)\times f_1(1+y) dy +o(\eps),\\
\label{d9} \tr^{(2)}(x_1,t_1,x_2,t_2)=\eps g^{(2)}(x_1,t_1,x_2,t_2)\\
+\displaystyle{\int_{-1}^\infty} \tr(x_1,t_1,x_2+\delta
y,t_2+\eps) \times f_2(1+y) dy +o(\eps). \nonumber
\end{gather}
Bayesian strategy prescribes to choose  $\ell$-th action if
$\tr^{(\ell)}(x_1,n_1,x_2,n_2)$ has smaller value. In case of a
draw the choice of the action is arbitrary. For the strategy which
at the start of the control equally applies both actions
$2n_0=2\eps_0 N$ times and then optimally controls, one derives
that Bayesian risk~\eqref{a3} is calculated by the formula
\begin{gather}\label{d10}
   R_N(\mu)=(D N)^{1/2} \Big( \eps_0 \iint_\Theta |v_2-v_1| \rho(v_1,v_2) dv_1 dv_2\\
   + \iint_{-\infty}^\infty \tr(\eps_0,x_1,\eps_0,x_2) dx_1
   dx_2\Big).\nonumber
\end{gather}
Formula~\eqref{d10} follows from~\eqref{c6} with the use
of~\eqref{d1} and~\eqref{d6}.

Finally, let's present a limiting description of \eqref{d9} by the
second order partial differential equation. It is sufficient to
consider the first equation of  \eqref{d9}. The estimates below
are carried out with accuracy to the terms of the order $\eps$.
First, we present $\tr(x_1+\delta y,t_1+\eps,x_2,t_2)$ as Taylor
series
\begin{gather}\label{d11}
\tr(x_1+\delta y,t_1+\eps,x_2,t_2)=\tr(x_1,t_1+\eps,x_2,t_2)\\
+\delta \tr'_{x_1}(x_1,t_1+\eps,x_2,t_2) y+0.5 \eps \tr''_{x_1
x_1}(x_1,t_1+\eps,x_2,t_2)y^2 +o(\eps). \nonumber
\end{gather}
Substituting \eqref{d11} into the first equation \eqref{d9} and
using \eqref{d5} we obtain
\begin{gather*}
\tr^{(1)}(x_1,t_1,x_2,t_2)=\eps g^{(1)}(x_1,t_1,x_2,t_2)+\tr(x_1,t_1+\eps,x_2,t_2)(1+\eps/t_1)\\
+\eps \tr'_{x_1}(x_1,t_1+\eps,x_2,t_2)\tx_1+0.5 \eps \tr''_{x_1
x_1}(x_1,t_1+\eps,x_2,t_2)+o(\eps),
\end{gather*}
and, hence, in the limiting case as $\eps \to 0$
\begin{gather}\label{d12}
g^{(1)}(x_1,t_1,x_2,t_2)+\tr/t_1+ \tr'_{t_1} + \tr'_{x_1}
\tx_1+0.5 \tr''_{x_1 x_1}=0,
\end{gather}
where $\tr=\tr(x_1,t_1,x_2,t_2)$. Similarly,
\begin{gather}\label{d13}
g^{(2)}(x_1,t_1,x_2,t_2)+\tr/t_2+ \tr'_{t_2} + \tr'_{x_2}
\tx_2+0.5 \tr''_{x_2 x_2}=0.
\end{gather}

Equations~\eqref{d12}, \eqref{d13} must be complemented by
equation \eqref{d7}, which is now written as
\begin{gather} \label{d14}
\min_{\ell=1,2}
(\tr(x_1,t_1,x_2,t_2)-\tr^{(\ell)}(x_1,t_1,x_2,t_2))=0
\end{gather}

From~\eqref{d12}--\eqref{d14} one derives
\begin{gather}\label{d15}
\min_{\ell=1,2}
\left(g^{(\ell)}(x_1,t_1,x_2,t_2)+\tr/t_\ell+\tr'_{t_\ell} +
\tr'_{x_\ell} \tx_\ell+0.5 \tr''_{x_\ell x_\ell}\right)=0.
\end{gather}
Initial  conditions are the following
\begin{gather}\label{d16}
r(x_1,t_1,x_2,t_2)=0,
\end{gather}
if $t_1+t_2=1$. Bayesian strategy prescribes to choose  $\ell$-th
action if $\ell$-th term in the left-hand side of \eqref{d15} has
smaller value. In case of a draw the choice of the action is
arbitrary.

\section{Comparison with Gaussian two-armed bandit}\label{Gaussian}

Gaussian two-armed bandit is characterized by incomes $\xi_n$,
$n=1,2,\dots,N$, which values depend only on currently chosen
actions $y_n \in \{1,2\}$ and are described by Gaussian (normal)
distribution density
\begin{gather}\label{e1}
f_D(x|m_\ell)= (2 \pi D)^{-1/2}
\exp\left(-(x-m_\ell)^2/(2D)\right),
\end{gather}
if $y_n=\ell$, $\ell=1,2$.  The variance $D$ is assumed to be
known and expectations $m_1, m_2$ are unknown. So, a vector
parameter $\theta=(m_1,m_2)$ describes Gaussian two-armed bandit.
We assume that the set of admissible values of parameters
$\Theta=\{\theta\}$ is a priori known.

A control strategy $\sigma$ generally assigns a random choice of
the action at the point of time $n+1$ depending on currently
observed history $(X_1,n_1,X_2,n_2)$ where $n_1, n_2$
($n_1+n_2=n$) are  cumulative times of both actions applications
and $X_1, X_2$ are corresponding cumulative incomes.

Again, one can assign a prior distribution density $\mu(m_1,m_2)$
and define a regret $L_N(\sigma,\theta)$  and Bayesian risk
$R_N(\mu)$ just like in \eqref{a2} and \eqref{a3}. To determine
Bayesian risk in the domain of ``close distributions'' one should
solve the following integro-difference equation (see, e.g.
\cite{Koln18}):
\begin{gather}\label{e2}
\tr(x_1,t_1,x_2,t_2)= \min (\tr^{(1)}(x_1,t_1,x_2,t_2),
\tr^{(2)}(x_1,t_1,x_2,t_2)),
\end{gather}
where
\begin{gather}\label{e3}
\tr^{(1)}(x_1,t_1,x_2,t_2)= \tr^{(2)}(x_1,t_1,x_2,t_2)=0
\end{gather}
if $t_1+t_2=1$ and then
\begin{gather} \nonumber
\tr^{(1)}(x_1,t_1,x_2,t_2)=\eps g^{(1)}(x_1,t_1,x_2,t_2)\\
\nonumber + (t_1+\eps) \displaystyle{ \int_{-\infty}^\infty}
\tr(x_1+\delta
y,t_1+\eps,x_2,t_2)\times  f_{t_1(t_1+\eps)}(\eps x_1-t_1 y) dy +o(\eps),\\
\label{e4} \tr^{(2)}(x_1,t_1,x_2,t_2)=\eps g^{(2)}(x_1,t_1,x_2,t_2)\\
+ (t_2+\eps) \displaystyle{\int_{-\infty}^\infty}
\tr(x_1,t_1,x_2+\delta y,t_2+\eps) \times f_{t_2(t_2+\eps)}(\eps
x_2-t_2 y)dy +o(\eps). \nonumber
\end{gather}
Here $f_D(x)=f_D(x|0)$ and $x_1,x_2,t_1,t_2,\eps,\delta$ are
defined in Section~\ref{lim}. Bayesian strategy prescribes to
choose the $\ell$-th action if $\tr^{(\ell)}(x_1,n_1,x_2,n_2)$ has
smaller value. In case of a draw the choice of the action is
arbitrary. For the strategy which at the start of the control
equally applies both actions $2n_0=2\eps_0 N$ times and then
optimally controls, one derives that Bayesian risk~\eqref{a3} is
calculated by the formula
\begin{gather}\label{e5}
   R_N(\mu)=(D N)^{1/2} \Big( \eps_0 \iint_\Theta |v_2-v_1| \rho(v_1,v_2) dv_1 dv_2\\
   + \iint_{-\infty}^\infty \tr(\eps_0,x_1,\eps_0,x_2) dx_1
   dx_2\Big).\nonumber
\end{gather}

Here $\rho(v_1,v_2)$ is defined in Section~\ref{lim}. In the
limiting case as $\eps \to 0$, one can verify that
integro-difference equation \eqref{e4} results in the second order
partial differential equation~\eqref{d15} with initial
conditions~\eqref{d16}. Recall that Gaussian two-armed bandit
describes the batch processing~\cite{Koln18}. So, this means that
in case of exponential two-armed bandit batch processing does not
enlarge Bayesian risk in comparison with one-by-one optimal
processing asymptotically as $N \to \infty$.



\vspace{6pt}

\begin{thebibliography}{9}

\bibitem{BF}
\textsc{Berry, D. A.} and \textsc{Fristedt, B.} (1985).
\textit{Bandit Problems: Sequential Allocation of Experiments},
Chapman~\&\ Hall, London.

\bibitem{PS}
\textsc{Presman, E. L.} and \textsc{Sonin, I. M.} (1990).
\textit{Sequential Control with Incomplete Information: Bayesian
Approach}, Academic Press, New York.

\bibitem{Sragovich}
\textsc{Sragovich, V. G.} (2006). \textit{Mathematical Theory of
Adaptive Control}, World Sci., Singapore.

\bibitem{CBL}
\textsc{Cesa-Bianchi, N.} and \textsc{Lugosi. G.} (2006)
\textit{Prediction, Learning, and Games}, Cambridge Univ. Press,
Cambridge.

\bibitem{Koln18}
\textsc{Kolnogorov, A. V.} (2018). Gaussian Two-Armed Bandit and
Optimization of Batch Data Processing. \textit{Problems of
Information Transmission} \textbf{54} 84--100.


\end{thebibliography}
\end{document}